\documentclass[11pt]{amsart}
\usepackage{amssymb,latexsym,amsmath,graphicx,graphics,epic,eepic}

\theoremstyle{plain}
\newtheorem{theorem}{Theorem}
\numberwithin{theorem}{section}
\newtheorem{lemma}[theorem]{Lemma}

\theoremstyle{definition}

\newtheorem{remark}[theorem]{Remark}

\begin{document}

\title{Tight Contact Small Seifert Spaces with $e_0\neq0,-1,-2$}

\author{Hao Wu}

\address{Department of mathematics\\
         MIT, Room 2-487\\
         77 Massachusetts Avenue\\
         Cambridge, MA 02139\\
         USA}

\email{haowu@math.mit.edu}

\begin{abstract}
We classify up to isotopy the tight contact structures on small
Seifert spaces with $e_0\neq0,-1,-2$.
\end{abstract}

\maketitle

\section{Introduction and Statements of Results}

A contact structure $\xi$ on an oriented $3$-manifold $M$ is a
nowhere integrable tangent plane distribution, i.e., near any
point of $M$, $\xi$ is defined locally by a $1$-form $\alpha$,
s.t., $\alpha\wedge d\alpha\neq0$. Note that the orientation of
$M$ given by $\alpha\wedge d\alpha$ depends only on $\xi$, not on
the choice of $\alpha$. $\xi$ is said to be positive if this
orientation agrees with the native orientation of $M$, and
negative if not. A contact structure $\xi$ is said to be
co-orientable if $\xi$ is defined globally by a $1$-form $\alpha$.
Clearly, an co-orientable contact structure is orientable as a
plane distribution, and a choice of $\alpha$ determines an
orientation of $\xi$. Unless otherwise specified, all contact
structures in this paper will be co-oriented and positive, i.e.,
with a prescribed defining form $\alpha$ such that $\alpha\wedge
d\alpha>0$. A curve in $M$ is said to be Legendrian if it is
tangent to $\xi$ everywhere. $\xi$ is said to be overtwisted if
there is an embedded disk $D$ in $M$ such that $\partial D$ is
Legendrian, but $D$ is transversal to $\xi$ along $\partial D$. A
contact structure that is not overtwisted is called tight.
Overtwisted contact structures are classified up to isotopy by
Eliashberg in \cite{E1}. Classifying tight contact structures up
to isotopy is much more difficult. Such classifications are only
known for very limited classes of $3$-manifolds. (See, e.g.,
\cite{E2}, \cite{E3}, \cite{EH}, \cite{Gh}, \cite{GS}, \cite{H1},
\cite{H2}, \cite{K}.)

For a small Seifert space $M=M(r_1,r_2,r_3)$, define $e_0(M) =
\lfloor r_1 \rfloor + \lfloor r_2 \rfloor + \lfloor r_3 \rfloor$,
where $\lfloor x \rfloor$ is the greatest integer not greater than
$x$. $e_0(M)$ is an invariant of $M$, i.e., it does not depend on
the choice of the representatives $(r_1,r_2,r_3)$.

For a rational number $r>0$, there is a unique way to expand $-r$
into a continued fraction
\begin{equation}\label{expansion}
-r ~=~ a_0 - \frac{1}{a_1 - \frac{1}{a_2 - \cdots\frac{1}{a_{m-1}
- \frac{1}{a_m}}}},
\end{equation}
where all $a_j$'s are integers, $a_0=-(\lfloor r \rfloor
+1)\leq-1$, and $a_j\leq-2$ for $j\geq1$. We denote by
$<a_0,a_1,\cdots,a_m>$ the right hand side of equation
\eqref{expansion}.

The following are our main results.

\begin{theorem}\label{class-3}
Let $M=M(-\frac{q_1}{p_1},-\frac{q_2}{p_2},-\frac{q_3}{p_3})$ be a
small Seifert space, where $p_i$ and $q_i$ are integers, s.t.,
$p_i\geq2$, $q_i\geq1$ and $g.c.d.(p_i,q_i)=1$. Assume that, for
$i=1,2,3$,
$-\frac{q_i}{p_i}=<a^{(i)}_0,a^{(i)}_1,\cdots,a^{(i)}_{m_i}>$,
where all $a^{(i)}_j$'s are integers,
$a^{(i)}_0=-(\lfloor\frac{q_i}{p_i}\rfloor+1)\leq-1$, and
$a^{(i)}_j\leq-2$ for $j\geq1$. Then, up to isotopy, there are
exactly $|(e_0(M)+1)\prod^3_{i=1}\prod^{m_i}_{j=1}(a^{(i)}_j+1)|$
tight contact structures on $M$. All these tight contact
structures are constructed by Legendrian surgeries of
$(S^3,\xi_{st})$, and are therefore holomorphically fillable.
\end{theorem}

\begin{theorem}\label{class+1}
Let $M=M(\frac{q_1}{p_1},\frac{q_2}{p_2},e_0+\frac{q_3}{p_3})$ be
a small Seifert space, where $p_i$, $q_i$ and $e_0$ are positive
integers, s.t., $p_i>q_i$ and $g.c.d.(p_i,q_i)=1$. Assume that,
for $i=1,2,3$,
$-\frac{p_i}{q_i}=<b^{(i)}_0,b^{(i)}_1,\cdots,b^{(i)}_{l_i}>$,
where all $b^{(i)}_j$'s are integers less than or equal to $-2$.
Then, up to isotopy, there are exactly
$|\prod^3_{i=1}b^{(i)}_0\prod^{l_i}_{j=1}(b^{(i)}_j+1)|$ tight
contact structures on $M$. All these tight contact structures are
constructed by Legendrian surgeries of $(S^3,\xi_{st})$, and are
therefore holomorphically fillable.
\end{theorem}

Clearly, Theorems \ref{class-3} and \ref{class+1} give complete
classifications of tight contact structures on small Seifert
spaces with $e_0\neq 0,-1,-2$.

\vspace{.4cm}

\section{Continued Fractions}

In this section, we establish some properties of continued
fractions, which will be used in the proofs of Theorems
\ref{class-3} and \ref{class+1}.

\begin{lemma}\label{CF}
Let $a_0,~a_1,\cdots,~a_m$ be real numbers such that $a_0\leq-1$,
and $a_j\leq-2$ for $1\leq j \leq m$. Define $\{p_j\}$ and
$\{q_j\}$ by
\[
\left\{%
\begin{array}{ll}
    p_j=-a_jp_{j-1}-p_{j-2}, ~j=0,1,\cdots,m, \\
    p_{-2}=-1, ~p_{-1}=0,  \\
\end{array}%
\right.
\]
\[
\left\{%
\begin{array}{ll}
    q_j=-a_jq_{j-1}-q_{j-2}, ~j=0,1,\cdots,m, \\
    q_{-2}=0, ~q_{-1}=1. \\
\end{array}%
\right.
\]
Then, for $1 \leq j \leq m$, we have
\begin{enumerate}
    \item $-\frac{q_j}{p_j}=<a_0,a_1,\cdots,a_j>$,
    \item $p_j \geq p_{j-1}>0$, $q_j \geq q_{j-1}>0$,
    \item $p_jq_{j-1}-p_{j-1}q_j=1$,
    \item
    $-\frac{q_j+(a_0+1)p_j}{q_{j-1}+(a_0+1)p_{j-1}}=<a_j,a_{j-1},\cdots,a_2,a_1+1>$.
\end{enumerate}
\end{lemma}

\begin{proof}
By the definitions of $\{p_j\}$ and $\{q_j\}$, we have $p_0=1$,
$q_0=-a_0$, $p_1=-a_1$, and $q_1=a_0a_1-1$. Then it's easy to
check that the lemma is true for $j=1$. Assume that the lemma is
true for $j-1\geq1$. Then,
\begin{eqnarray*}
<a_0,a_1,\cdots,a_j> & = & <a_0,a_1,\cdots,a_{j-1}-\frac{1}{a_j}>
\\
& = &
-~\frac{-(a_{j-1}-\frac{1}{a_j})q_{j-2}-q_{j-3}}{-(a_{j-1}-\frac{1}{a_j})p_{j-2}-p_{j-3}}
\\
& = & -~
\frac{(a_ja_{j-1}-1)q_{j-2}+a_jq_{j-3}}{(a_ja_{j-1}-1)p_{j-2}+a_jp_{j-3}}
\\
& = & -~
\frac{a_j(a_{j-1}q_{j-2}+q_{j-3})-q_{j-2}}{a_j(a_{j-1}p_{j-2}+p_{j-3})-p_{j-2}}
\\
& = &  -~\frac{-a_jq_{j-1}-q_{j-2}}{-a_jq_{j-1}-q_{j-2}} \\
& = & -~\frac{q_j}{p_j}.
\end{eqnarray*}

Also, since $q_{j-1}\geq q_{j-2}>0$ and $-a_j\geq2$, we have
$q_j=-a_jq_{j-1}-q_{j-2}\geq2q_{j-1}-q_{j-2}\geq q_{j-1}>0$, and,
similarly, $p_j\geq p_{j-1}>0$.

Furthermore, by definitions of $\{p_j\}$ and $\{q_j\}$,
\begin{eqnarray*}
p_jq_{j-1}-p_{j-1}q_j & = &
(-a_jp_{j-1}-p_{j-2})q_{j-1}-p_{j-1}(-a_jq_{j-1}-q_{j-2}) \\
& = & p_{j-1}q_{j-2}-p_{j-2}q_{j-1} \\
& = & 1.
\end{eqnarray*}

Finally,
\begin{eqnarray*}
-~\frac{q_j+(a_0+1)p_j}{q_{j-1}+(a_0+1)p_{j-1}} & = &
-~\frac{(-a_jq_{j-1}-q_{j-2})+(a_0+1)(-a_jp_{j-1}-p_{j-2})}{q_{j-1}+(a_0+1)p_{j-1}}
\\
& = &
\frac{a_j(q_{j-1}+(a_0+1)p_{j-1})+(q_{j-2}+(a_0+1)p_{j-2})}{q_{j-1}+(a_0+1)p_{j-1}}
\\
& = & a_j - \frac{1}{<a_{j-1},\cdots,a_2,a_1+1>} \\
& = & <a_j,a_{j-1},\cdots,a_2,a_1+1>.
\end{eqnarray*}

This shows that the lemma is also true for $j$.
\end{proof}

\begin{remark}\label{integer}
In the proof of Theorem \ref{class-3} and \ref{class+1}, all the
$a_j$'s will be integers, and so will the corresponding $p_j$'s
and $q_j$'s be. Then, property (3) in Lemma \ref{CF} implies that
$g.c.d.(p_j,q_j)=1$.
\end{remark}

\vspace{.4cm}

\section{The $e_0\leq-3$ Case}

In the rest of this paper, we let $\Sigma$ be a three hole sphere,
and $-\partial\Sigma\times S^1=T_1+T_2+T_3$, where the "$-$" sign
means reversing the orientation. We identify $T_i$ to
$\mathbb{R}^2/\mathbb{Z}^2$ by identifying the corresponding
component of $-\partial\Sigma\times\{\text{pt}\}$ to $(1,0)^T$,
and $\{\text{pt}\}\times S^1$ to $(0,1)^T$. Also, for $i=1,2,3$,
let $V_i=D^2\times S^1$, and identify $\partial V_i$ with
$\mathbb{R}^2/\mathbb{Z}^2$ by identifying a meridian $\partial
D^2 \times \{\text{pt}\}$ with $(1,0)^T$ and a longitude
$\{\text{pt}\}\times S^1$ with $(0,1)^T$.

Following Honda, we call a convex torus minimal if it has only two
dividing curves.

\vspace{.4cm}

\noindent\textit{Proof of Theorem \ref{class-3}.} Define
$\{p^{(i)}_j\}$ and $\{q^{(i)}_j\}$ by
\[
\left\{%
\begin{array}{ll}
    p^{(i)}_j=-a^{(i)}_jp^{(i)}_{j-1}-p^{(i)}_{j-2}, ~j=0,1,\cdots,m_i, \\
    p^{(i)}_{-2}=-1, ~p^{(i)}_{-1}=0,  \\
\end{array}%
\right.
\]
\[
\left\{%
\begin{array}{ll}
    q^{(i)}_j=-a^{(i)}_jq^{(i)}_{j-1}-q^{(i)}_{j-2}, ~j=0,1,\cdots,m_i, \\
    q^{(i)}_{-2}=0, ~q^{(i)}_{-1}=1. \\
\end{array}%
\right.
\]
By Lemma \ref{CF} and Remark \ref{integer}, we have
$p_i=p^{(i)}_{m_i}$ and $q_i=q^{(i)}_{m_i}$. Let
$u_i=p^{(i)}_{m_i-1}$ and $v_i=q^{(i)}_{m_i-1}$. Then $p_i \geq
u_i>0$, $q_i \geq v_i>0$, and $p_iv_i-q_iu_i=1$.

Define an orientation
preserving diffeomorphism $\varphi_i:\partial V_i\rightarrow T_i$
by
\[
\varphi_i ~=~
\left(%
\begin{array}{cc}
  p_i & u_i \\
  q_i & v_i \\
\end{array}%
\right).
\]
Then
\[
M=M(-\frac{q_1}{p_1},-\frac{q_2}{p_2},-\frac{q_3}{p_3})\cong(\Sigma
\times S^1)\cup_{(\varphi_1\cup\varphi_2\cup\varphi_3)}(V_1\cup
V_2\cup V_3).
\]

Let $\xi$ be a tight contact structure on $M$. We first isotope
$\xi$ to make each $V_i$ a standard neighborhood of a Legendrian
circle $L_i$ isotopic to the $i$-th singular fiber with twisting
number $t_i<-2$, i.e., $\partial V_i$ is convex with two dividing
curves each of which has slope $\frac{1}{t_i}$ when measured in
the coordinates of $\partial V_i$ given above. Let $s_i$ be the
slope of the dividing curves of $T_i=\partial V_i$ measured in the
coordinates of $T_i$. Then we have that
\[
s_i ~=~ \frac{t_iq_i+v_i}{t_ip_i+u_i} ~=~
\frac{q_i}{p_i}+\frac{1}{p_i(t_ip_i+u_i)}.
\]
The fact $t_i<-2$ implies that
$\lfloor\frac{q_i}{p_i}\rfloor<s_i<\frac{q_i}{p_i}$.

After a possible slight isotopy supported in a neighborhood of
$T_i=\partial V_i$, we assume that $T_i$ has Legendrian ruling of
slope $\infty$ when measured in the coordinates of $T_i$. For each
$i$, pick a Legendrian ruling $L_i$ on $T_i$. Choose a convex
vertical annulus $A \subset \Sigma \times S^1$, such that
$\partial A = L_1 \cup L_2$, and the interior of $A$ is contained
in the interior of $\Sigma \times S^1$. By Theorem 1.4 of
\cite{Wu}, $\xi$ does not admit Legendrian vertical circles with
twisting number $0$. So there must be dividing curves of $A$ that
connect the two boundary components of $A$. We isotope $T_1$ and
$T_2$ by adding to them the bypasses corresponding to the
$\partial$-parallel dividing curves of $A$. Since bypass adding is
done in a small neighborhood of the bypass and the original
surface, we can keep $V_i$'s disjoint during this process. Also
$T_i$ remains minimal after each bypass adding. After we depleted
all the $\partial$-parallel dividing curves of $A$, each of the
remaining dividing curves connects the two boundary components of
$A$. So the slopes of the dividing curves of $T_1$ and $T_2$ after
the isotopy are $s'_1=\frac{k_1}{k}$ and $s'_2=\frac{k_2}{k}$,
where $k\geq1$ and $g.c.d.(k,k_i)=1$ for $i=1,2$. Since
$\lfloor\frac{q_i}{p_i}\rfloor<s_i$, We have that, for $i=1,2$,
$s'_i\geq\lfloor\frac{q_i}{p_i}\rfloor\geq0$, and, hence
$k_i\geq0$. This is because that, by Lemma 3.15 of \cite{H1},
$s'_i<\lfloor\frac{q_i}{p_i}\rfloor$ implies $s'_i=\infty$ which
contradicts Theorem 1.4 of \cite{Wu}. Now, cut $M$ open along $A
\cup T_1 \cup T_2$ and round the edges. We get a convex torus
isotopic to $T_3$ with two dividing curves of slope
$-\frac{k_1+k_2+1}{k}$ when measure in the coordinates of $T_3$.
When measured in the coordinates of $\partial V_3$, these dividing
curves have slope
$-\frac{kq_3+(k_1+k_2+1)p_3}{kv_3+(k_1+k_2+1)u_3}$. It's easy to
check that
$-\frac{kq_3+(k_1+k_2+1)p_3}{kv_3+(k_1+k_2+1)u_3}<-\frac{q_3}{v_3}$.
So, by Theorem 4.16 of \cite{H1}, we can isotope $\partial V_3$ so
that it has two dividing curves of slope $-\frac{q_3}{v_3}$.
Measured in the coordinates of $T_3$, the slope is $0$. This
implies that the maximal twisting number of a Legendrian vertical
circle is $-1$.

After an isotopy of $\xi$, we can find a Legendrian vertical
circle $L$ in the interior of $\Sigma \times S^1$ with twisting
number $-1$, and, again, make each $V_i$ a standard neighborhood
of a Legendrian circle $L_i$ isotopic to the $i$-th singular fiber
with twisting number $t_i<-2$. As before, we can assume that $T_i$
has Legendrian ruling of slope $\infty$ when measured in the
coordinates of $T_i$. Let $L_i$ be a Legendrian ruling of $T_i$.
For each $i$, we choose a convex vertical annulus $A_i \subset
\Sigma \times S^1$, s.t., $\partial A_i = L \cup L_i$, the
interior of $A_i$ is contained in the interior of $\Sigma \times
S^1$, and $A_i \cap A_j=L$ when $i\neq j$. $A_i$ has no
$\partial$-parallel dividing curves on the $L$ side since $t(L)$
is maximal. So the dividing set of $A_i$ consists of two curves
connecting $L$ to $L_i$ and possibly some $\partial$-parallel
curves on the $L_i$ side. We now isotope $T_i$ by adding to it the
bypasses corresponding to these $\partial$-parallel dividing
curves, and keep $V_i$'s disjoint in this process. After this
isotopy, we get a convex decomposition
\[
M=M(-\frac{q_1}{p_1},-\frac{q_2}{p_2},-\frac{q_3}{p_3})\cong(\Sigma
\times S^1)\cup_{(\varphi_1\cup\varphi_2\cup\varphi_3)}(V_1\cup
V_2\cup V_3)
\]
of $M$, where each $T_i$ has two dividing curves of slope
$\lfloor\frac{q_i}{p_i}\rfloor$ when measured in the coordinate of
$T_i$. When measured in coordinates of $\partial V_i$, the slope
of the dividing curves becomes
$-\frac{q_i-\lfloor\frac{q_i}{p_i}\rfloor
p_i}{v_i-\lfloor\frac{q_i}{p_i}\rfloor u_i}=
-\frac{q_i+(a^{(i)}_0+1)p_i}{v_i+(a^{(i)}_0+1)u_i}$.

By part (4) of Lemma 5.1 of \cite{H2}, there are exactly $2+
\lfloor\frac{q_1}{p_1}\rfloor + \lfloor\frac{q_2}{p_2}\rfloor +
\lfloor\frac{q_3}{p_3}\rfloor = |e_0(M)+1|$ tight contact
structures on $\Sigma \times S^1$ satisfying the boundary
condition and admitting no Legendrian vertical circles with
twisting number $0$. By Theorem 2.3 of \cite{H1} and part (4) of
Lemma \ref{CF}, there are exactly
$|\prod^{m_i}_{j=1}(a^{(i)}_j+1)|$ tight contact structures on
$V_i$ satisfying the boundary condition. Thus, up to isotopy,
there are at most
$|(e_0(M)+1)\prod^3_{i=1}\prod^{m_i}_{j=1}(a^{(i)}_j+1)|$ tight
contact structures on $M$.

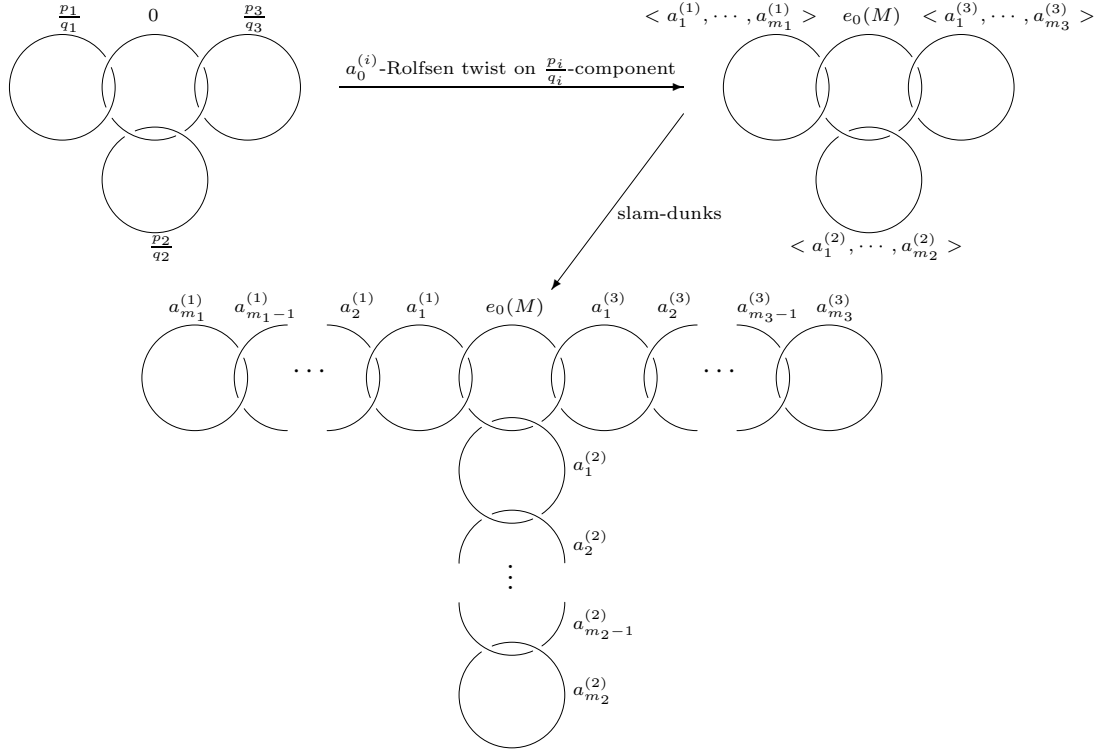
\begin{figure}[h]

\setlength{\unitlength}{1pt}

\begin{picture}(420,290)(-210,-250)


\put(-137.5,25){\tiny{$0$}}

\put(-135,0){\arc{40}{9.04}{11.95}}

\put(-135,0){\arc{40}{5.9}{7.25}}

\put(-135,0){\arc{40}{7.45}{8.8}}

\put(-172.5,25){\tiny{$\frac{p_1}{q_1}$}}

\put(-170,0){\arc{40}{5.9}{11.95}}

\put(-102.5,25){\tiny{$\frac{p_3}{q_3}$}}

\put(-100,0){\arc{40}{2.76}{8.8}}

\put(-137.5,-62){\tiny{$\frac{p_2}{q_2}$}}

\put(-135,-35){\arc{40}{4.33}{10.38}}


\put(-63,5){\tiny{$a^{(i)}_0$-Rolfsen twist on
$\frac{p_i}{q_i}$-component}}

\put(-65,0){\vector(1,0){130}}


\put(125,25){\tiny{$e_0(M)$}}

\put(135,0){\arc{40}{9.04}{11.95}}

\put(135,0){\arc{40}{5.9}{7.25}}

\put(135,0){\arc{40}{7.45}{8.8}}

\put(50,25){\tiny{$<a^{(1)}_1,\cdots,a^{(1)}_{m_1}>$}}

\put(100,0){\arc{40}{5.9}{11.95}}

\put(155,25){\tiny{$<a^{(3)}_1,\cdots,a^{(3)}_{m_3}>$}}

\put(170,0){\arc{40}{2.76}{8.8}}

\put(105,-62){\tiny{$<a^{(2)}_1,\cdots,a^{(2)}_{m_2}>$}}

\put(135,-35){\arc{40}{4.33}{10.38}}


\put(65,-10){\vector(-3,-4){50}}

\put(40,-50){\tiny{slam-dunks}}



\put(-10,-85){\tiny{$e_0(M)$}}

\put(0,-110){\arc{40}{9.04}{11.95}}

\put(0,-110){\arc{40}{5.9}{7.25}}

\put(0,-110){\arc{40}{7.45}{8.8}}


\put(-40,-85){\tiny{$a^{(1)}_1$}}

\put(-35,-110){\arc{40}{9.04}{11.95}}

\put(-35,-110){\arc{40}{5.9}{8.8}}

\put(-65,-85){\tiny{$a^{(1)}_2$}}

\put(-70,-110){\arc{40}{11}{11.95}}

\put(-70,-110){\arc{40}{5.9}{7.85}}

\put(-83,-110){$\cdots$}

\put(-105,-85){\tiny{$a^{(1)}_{m_1-1}$}}

\put(-85,-110){\arc{40}{9.04}{11}}

\put(-85,-110){\arc{40}{7.85}{8.8}}

\put(-130,-85){\tiny{$a^{(1)}_{m_1}$}}

\put(-120,-110){\arc{40}{5.9}{11.95}}


\put(30,-85){\tiny{$a^{(3)}_1$}}

\put(35,-110){\arc{40}{9.04}{11.95}}

\put(35,-110){\arc{40}{5.9}{8.8}}

\put(55,-85){\tiny{$a^{(3)}_2$}}

\put(70,-110){\arc{40}{9.04}{11}}

\put(70,-110){\arc{40}{7.85}{8.8}}

\put(72,-110){$\cdots$}

\put(85,-85){\tiny{$a^{(3)}_{m_3-1}$}}

\put(85,-110){\arc{40}{11}{11.95}}

\put(85,-110){\arc{40}{5.9}{7.85}}

\put(115,-85){\tiny{$a^{(3)}_{m_3}$}}

\put(120,-110){\arc{40}{9.04}{15.1}}


\put(23,-145){\tiny{$a^{(2)}_1$}}

\put(0,-145){\arc{40}{4.33}{7.25}}

\put(0,-145){\arc{40}{7.45}{10.38}}

\put(23,-175){\tiny{$a^{(2)}_2$}}

\put(0,-180){\arc{40}{4.33}{6.28}}

\put(0,-180){\arc{40}{9.42}{10.38}}

\put(-2,-190){$\vdots$}

\put(23,-205){\tiny{$a^{(2)}_{m_2-1}$}}

\put(0,-195){\arc{40}{7.45}{9.42}}

\put(0,-195){\arc{40}{6.28}{7.25}}

\put(23,-230){\tiny{$a^{(2)}_{m_2}$}}

\put(0,-230){\arc{40}{4.33}{10.38}}

\end{picture}

  \caption{Construction of the tight contact structures on $M$}\label{surgery}
\end{figure}

It remains to construct
$|(e_0(M)+1)\prod^3_{i=1}\prod^{m_i}_{j=1}(a^{(i)}_j+1)|$ tight
contact structures on $M$ by Legendrian surgeries of
$(S^3,\xi_{st})$. We begin with the standard surgery diagram of
$M=M(-\frac{q_1}{p_1},-\frac{q_2}{p_2},-\frac{q_3}{p_3})$. Then,
for each $i$, perform an $a^{(i)}_0$-Rolfsen twist on the
$\frac{p_i}{q_i}$-component. Since
$a^{(1)}_0+a^{(2)}_0+a^{(3)}_0=e_0(M)$ and
$\frac{p_i}{q_i+a^{(i)}_0p_i}=<a^{(i)}_1,\cdots,a^{(i)}_{m_i}>$,
the new surgery coefficients of the four components are $e_0(M)$,
$<a^{(1)}_1,\cdots,a^{(1)}_{m_1}>$,
$<a^{(2)}_1,\cdots,a^{(2)}_{m_2}>$, and
$<a^{(3)}_1,\cdots,a^{(3)}_{m_3}>$. Now, we perform (inverses of)
the slam-dunks corresponding to the three continued fractions
here, which lead us to the diagram at the bottom of Figure
\ref{surgery}. Since the maximal Thurston-Bennequin number of an
unknot in $(S^3,\xi_{st})$ is $-1$, there are
$|(e_0(M)+1)\prod^3_{i=1}\prod^{m_i}_{j=1}(a^{(i)}_j+1)|$ ways to
realize this diagram by Legendrian surgeries. According to
Proposition 2.3 of \cite{Go} and Theorem 1.2 of \cite{LM}, each of
these Legendrian surgeries gives a non-isotopic tight contact
structure on $M$. Thus, up to isotopy, there are exactly
$|(e_0(M)+1)\prod^3_{i=1}\prod^{m_i}_{j=1}(a^{(i)}_j+1)|$ tight
contact structures on $M$. \hfill $\Box$

\vspace{.4cm}

\section{The $e_0\geq1$ Case}

Following \cite{Wu}, we call a tight contact structure on
$\Sigma\times S^1$ appropriate is there are no embedded thickened
tori in $\Sigma\times S^1$ with $I$-twisting $\geq\pi$ and
parallel to one of the boundary components.

The following lemma is a reformulation of parts (1), (2) and (3)
of Lemma 5.1 of \cite{H2}.

\begin{lemma}\label{center}
Let $\xi$ be an appropriate contact structure on $\Sigma\times
S^1$ with minimal convex boundary that admits a vertical
Legendrian circle with twisting number $0$. Assume that dividing
curves of $T_1$, $T_2$ and $T_3$ are of slopes $-1$, $-1$, $-n$,
respectively, where $n$ is an integer greater than $1$. Then there
is a factorization $\Sigma\times S^1=L_1 \cup L_2 \cup L_3
\cup(\Sigma'\times S^1)$, where $L_i$'s are embedded thickened
tori with minimal twisting and minimal convex boundary $\partial
L_i=T'_i-T_i$, s.t., dividing curves of $T'_i$ have slope
$\infty$. The appropriate contact structure $\xi$ is uniquely
determined by the signs of the basic slices $L_1$, $L_2$ and
$L_3$. The sign convention here is given by associating $(0,1)^T$
to $T'_i$.
\end{lemma}

\begin{proof}
We only prove the last sentence. The rest is just part (1) of
Lemma 5.1 of \cite{H2}. Let $\Sigma_0$ be a properly embedded
three hole sphere in $\Sigma\times S^1$ isotopic to
$\Sigma\times\{pt\}$, and $\Sigma'_0=\Sigma_0\cap(\Sigma'\times
S^1)$. We isotope $\Sigma_0$ so that $\Sigma_0$ and $\Sigma'_0$
are convex with Legendrian boundaries that intersect the dividing
curves of $\partial\Sigma\times S^1$ and $\partial\Sigma'\times
S^1$ efficiently. Then each component of $\partial\Sigma'_0$
intersects the dividing curves of $\Sigma'_0$ twice. Since $\xi$
is appropriate, $\Sigma'_0$ has no $\partial$-parallel dividing
curves. This implies that, up to isotopy relative to boundary and
Dehn twists parallel to boundary components, there are
only two configurations of dividing curves on $\Sigma'_0$. (See Figure \ref{inner}.)
Thus, there are only two tight contact structure on $\Sigma'\times S^1$,
up to isotopy relative to boundary and full horizontal rotations
of each boundary component.

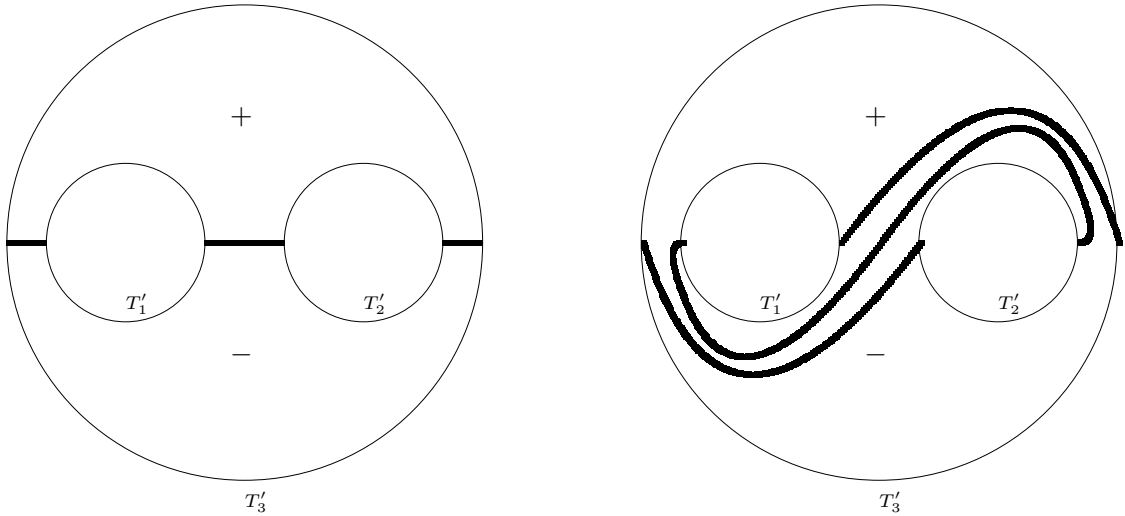
\begin{figure}[h]

\setlength{\unitlength}{1pt}

\begin{picture}(420,200)(-210,-100)

\put(-120,0){\circle{180}}

\put(-165,0){\circle{60}}

\put(-75,0){\circle{60}}

\put(-125.5,45){$+$}

\put(-125.5,-45){$-$}

\put(-120,-100){\tiny{$T'_3$}}

\put(-165,-25){\tiny{$T'_1$}}

\put(-75,-25){\tiny{$T'_2$}}

\put(120,0){\circle{180}}

\put(165,0){\circle{60}}

\put(75,0){\circle{60}}

\put(114.5,45){$+$}

\put(114.5,-45){$-$}

\put(120,-100){\tiny{$T'_3$}}

\put(165,-25){\tiny{$T'_2$}}

\put(75,-25){\tiny{$T'_1$}}

\linethickness{2pt}

\put(-210,0){\line(1,0){15}}

\put(-135,0){\line(1,0){30}}

\put(-45,0){\line(1,0){15}}

\qbezier(30,0)(60,-100)(135,0)

\qbezier(105,0)(180,100)(210,0)

\qbezier(195,0)(203,0)(195,20)

\qbezier(195,20)(175,75)(120,0)

\qbezier(45,0)(37,0)(45,-20)

\qbezier(45,-20)(65,-75)(120,0)

\end{picture}

\caption{Possible configurations of dividing curves on
$\Sigma'_0$}\label{inner}
\end{figure}

Let $A_i=\Sigma_0\cap L_i$. Then the dividing set of each of $A_1$
and $A_2$ consists of two arcs connecting the two boundary
components. And the dividing set of $A_3$ consists of two arcs
connecting the two boundary components and $n-1$
$\partial$-parallel arcs on the $T_3$ side. From the relative
Euler class of $\xi|_{L_3}$, one can see that the half discs
bounded by these $\partial$-parallel arcs must be pairwise
disjoint and of the sign opposite to that of $L_3$. By isotoping
$\Sigma_0$ relative to $\Sigma'_0$, we can freely choose the
holonomy of the non-$\partial$-parallel dividing curves of each
$A_i$. This implies that, up to isotopy relative to boundary,
there are only two possible configurations of dividing curves
on $\Sigma_0$ when the signs of $L_i$'s are given. (See Figure \ref{middle}.)

\begin{figure}[h]

\setlength{\unitlength}{1pt}

\begin{picture}(420,200)(-210,-100)

\put(-120,0){\circle{180}}

\put(-165,0){\circle{60}}

\put(-75,0){\circle{60}}

\put(-125.5,45){$+$}

\put(-125.5,-45){$-$}

\put(-180,55){$-$}

\put(-120,80){$-$}

\put(-120,-100){\tiny{$T_3$}}

\put(-165,-25){\tiny{$T_1$}}

\put(-75,-25){\tiny{$T_2$}}

\put(120,0){\circle{180}}

\put(165,0){\circle{60}}

\put(75,0){\circle{60}}

\put(114.5,45){$+$}

\put(114.5,-45){$-$}

\put(60,55){$-$}

\put(120,80){$-$}

\put(120,-100){\tiny{$T_3$}}

\put(165,-25){\tiny{$T_2$}}

\put(75,-25){\tiny{$T_1$}}

\linethickness{2pt}

\put(-210,0){\line(1,0){15}}

\put(-135,0){\line(1,0){30}}

\put(-45,0){\line(1,0){15}}

\qbezier(30,0)(60,-100)(135,0)

\qbezier(105,0)(180,100)(210,0)

\qbezier(195,0)(203,0)(195,20)

\qbezier(195,20)(175,75)(120,0)

\qbezier(45,0)(37,0)(45,-20)

\qbezier(45,-20)(65,-75)(120,0)

\qbezier(100,87)(120,55)(140,87)

\qbezier(43,46)(90,40)(83,81)

\qbezier(-140,87)(-120,55)(-100,87)

\qbezier(-197,46)(-150,40)(-157,81)

\end{picture}

\caption{Possible configurations of dividing curves on $\Sigma_0$.
Here, $n=3$, and the layer $L_3$ is positive}

\label{middle}

\end{figure}
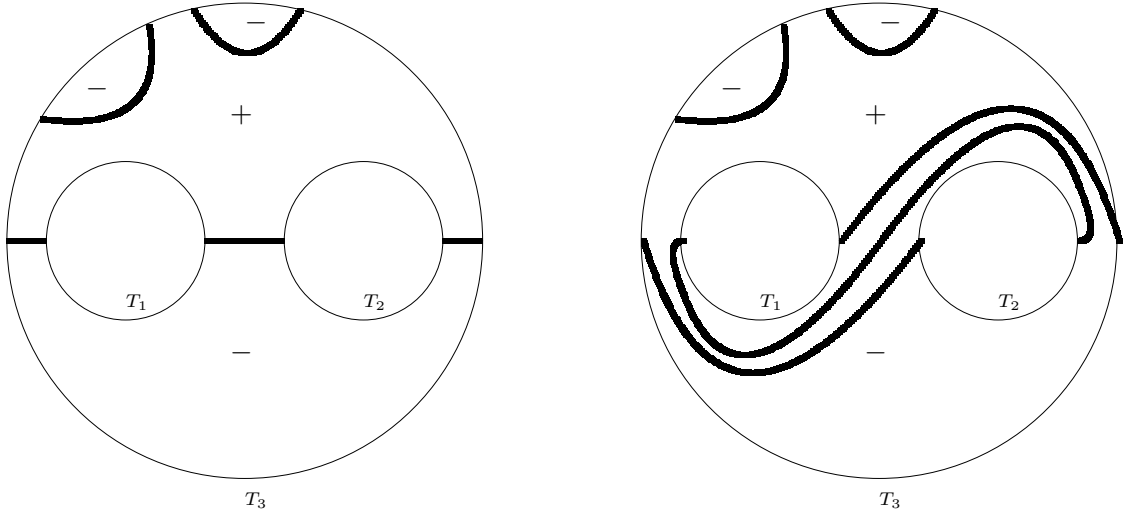

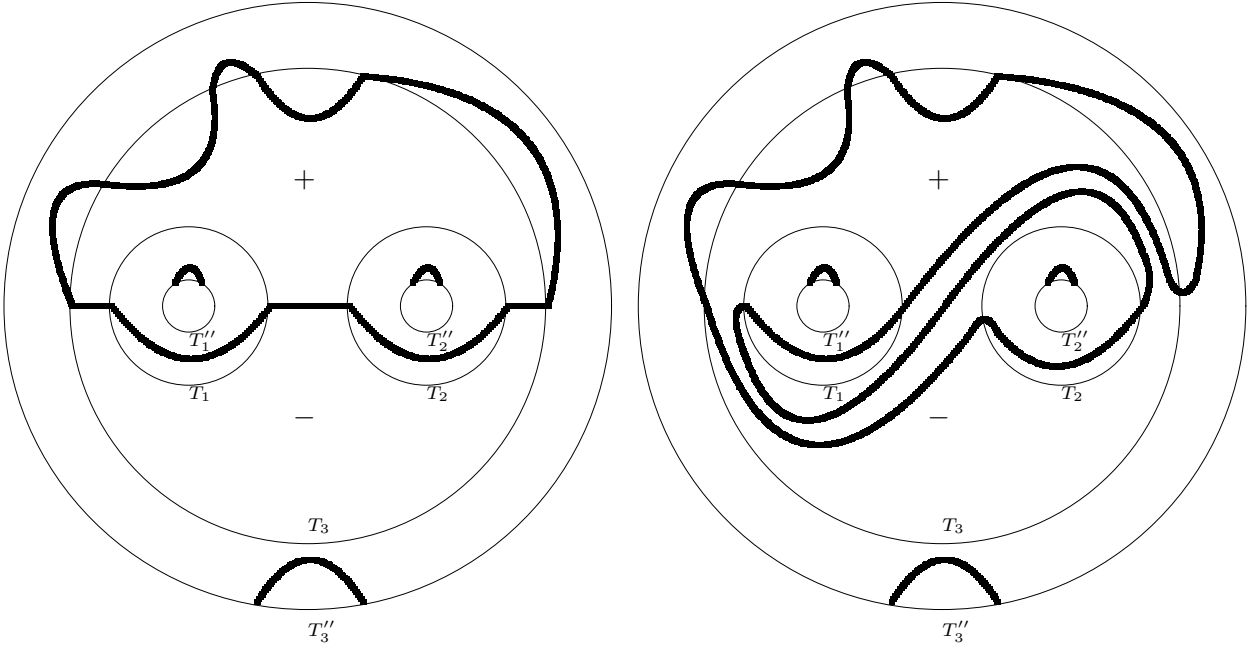
\begin{figure}[h]

\setlength{\unitlength}{1pt}

\begin{picture}(420,260)(-210,-130)

\put(-120,0){\circle{180}}

\put(-120,0){\circle{230}}

\put(-165,0){\circle{60}}

\put(-165,0){\circle{20}}

\put(-75,0){\circle{60}}

\put(-75,0){\circle{20}}

\put(-125.5,45){$+$}

\put(-125.5,-45){$-$}

\put(-120,-85){\tiny{$T_3$}}

\put(-165,-35){\tiny{$T_1$}}

\put(-75,-35){\tiny{$T_2$}}

\put(-120,-125){\tiny{$T''_3$}}

\put(-165,-15){\tiny{$T''_1$}}

\put(-75,-15){\tiny{$T''_2$}}

\put(120,0){\circle{180}}

\put(120,0){\circle{230}}

\put(165,0){\circle{60}}

\put(165,0){\circle{20}}

\put(75,0){\circle{60}}

\put(75,0){\circle{20}}

\put(114.5,45){$+$}

\put(114.5,-45){$-$}

\put(120,-85){\tiny{$T_3$}}

\put(165,-35){\tiny{$T_2$}}

\put(75,-35){\tiny{$T_1$}}

\put(120,-125){\tiny{$T''_3$}}

\put(165,-15){\tiny{$T''_2$}}

\put(75,-15){\tiny{$T''_1$}}

\linethickness{2pt}

\put(-210,0){\line(1,0){15}}

\put(-135,0){\line(1,0){30}}

\put(-45,0){\line(1,0){15}}

\qbezier(30,0)(60,-100)(130,-10)

\qbezier(105,0)(180,100)(205,10)

\qbezier(195,0)(200,10)(195,20)

\qbezier(195,20)(175,75)(120,0)

\qbezier(45,0)(37,0)(45,-20)

\qbezier(45,-20)(65,-75)(120,0)

\qbezier(100,87)(120,55)(140,87)

\qbezier(43,46)(90,40)(83,81)

\qbezier(-140,87)(-120,55)(-100,87)

\qbezier(-197,46)(-150,40)(-157,81)

\qbezier(83,81)(85,100)(100,87)

\qbezier(140,87)(230,80)(215,10)

\qbezier(205,10)(210,0)(215,10)

\qbezier(30,0)(10,50)(43,46)

\qbezier(45,0)(75,-40)(105,0)

\qbezier(140,-10)(165,-40)(195,0)

\qbezier(140,-10)(135,0)(130,-10)

\qbezier(100,-112)(120,-80)(140,-112)

\qbezier(69,9)(75,20)(79,9)

\qbezier(159,9)(165,20)(169,9)

\qbezier(-157,81)(-155,100)(-140,87)

\qbezier(-100,87)(-10,80)(-30,0)

\qbezier(-210,0)(-230,50)(-197,46)

\qbezier(-195,0)(-165,-40)(-135,0)

\qbezier(-105,0)(-75,-40)(-45,0)

\qbezier(-140,-112)(-120,-80)(-100,-112)

\qbezier(-171,9)(-165,20)(-161,9)

\qbezier(-81,9)(-75,20)(-71,9)

\end{picture}

\caption{After extending to $\Sigma''_0$, the two possible configurations become the same.
Here, $n=3$, and the signs of the layers $L_1$, $L_2$ and $L_3$ are $(-,-,+)$, respectively.}

\label{outer}

\end{figure}

When the signs of $L_i$'s are mixed, we can extend
$(\partial\Sigma\times S^1,\xi)$ to a universally tight contact
manifold $(\partial\Sigma''\times S^1,\xi'')$ by gluing to $T_i$ a
basic slice $L''_i$ of the same sign as $L_i$ for each $i$, where
$L''_i$ has minimal convex boundary $\partial L''_i=T_i-T''_i$,
and the dividing curves of $T''_i$ are vertical. Extend $\Sigma_0$
across $L''_i$ to $\Sigma''_0$ so that $\Sigma''_0$ is convex with
Legendrian boundary intersecting the dividing curves of $T''_i$
efficiently. For $i=1,2$, the dividing set of $\Sigma''_0\cap
L''_i$ consists of $1$ $\partial$-parallel arcs on each boundary
component. From the relative Euler class of $\xi''|_{L''_i}$, we
can see that the half discs on $\Sigma''_0\cap L''_i$ bounded by
these $\partial$-parallel arcs are of the same sign as the basic
slice $L_i$. The dividing set of $\Sigma''_0\cap L''_3$ consists
of $n$ $\partial$-parallel arcs on the $T_3$ side and $1$
$\partial$-parallel arcs on the $T''_3$ side. From the relative
Euler class of $\xi''|_{L''_3}$, we can see that the half discs on
$\Sigma''_0\cap L''_3$ bounded by these $\partial$-parallel arcs
are pairwise disjoint and of the same sign as the basic slice
$L_3$.

Now, one can see that, after the extension, the two
possible configurations of dividing curves on $\Sigma_0$ become
the same minimal configuration of dividing curves on $\Sigma''_0$. (See Figure \ref{outer}.)
By Lemma 5.1 of \cite{H2}, the two configurations correspond to
the same universally tight contact structure on $\Sigma\times
S^1$. This shows that, when the signs of $L_i$'s are mixed, $\xi$
is uniquely determined by the signs of $L_i$'s. When all the
$L_i$'s have the same sign, $\xi$ is virtually overtwisted, and
the isotopy type relative to boundary of such a
contact structure is determined by the action of the relative Euler
class on $\Sigma_0$, which is, in turn, determined by the sign of
$L_3$. Thus, when all the $L_i$'s have the same sign, this common
sign determines $\xi$.
\end{proof}

\vspace{.4cm}

\noindent\textit{Proof of Theorem \ref{class+1}.} Define
$\{p^{(i)}_j\}$ and $\{q^{(i)}_j\}$ by
\[
\left\{%
\begin{array}{ll}
    p^{(i)}_j=-b^{(i)}_jp^{(i)}_{j-1}-p^{(i)}_{j-2}, ~j=0,1,\cdots,l_i, \\
    p^{(i)}_{-2}=0, ~p^{(i)}_{-1}=1,  \\
\end{array}%
\right.
\]
\[
\left\{%
\begin{array}{ll}
    q^{(i)}_j=-b^{(i)}_jq^{(i)}_{j-1}-q^{(i)}_{j-2}, ~j=0,1,\cdots,l_i, \\
    q^{(i)}_{-2}=-1, ~q^{(i)}_{-1}=0. \\
\end{array}%
\right.
\]
By Lemma \ref{CF} and Remark \ref{integer}, we have
$p_i=p^{(i)}_{l_i}$ and $q_i=q^{(i)}_{l_i}$. Let
$u_i=-p^{(i)}_{l_i-1}$ and $v_i=-q^{(i)}_{l_i-1}$. Then
$p_iv_i-q_iu_i=1$.

Define an orientation preserving diffeomorphism
$\varphi_i:\partial V_i\rightarrow T_i$ by
\[
\varphi_i ~=~
\left\{%
\begin{array}{ll}
    \left(%
\begin{array}{cc}
  p_i & -u_i \\
  -q_i & v_i \\
\end{array}%
\right),
& i=1,2; \\
    \left(%
\begin{array}{cc}
  p_3 & -u_3 \\
  -q_3-e_0p_3 & v_3+e_0u_3 \\
\end{array}%
\right), & i=3. \\
\end{array}%
\right.
\]
Then
\[
M=M(\frac{q_1}{p_1},\frac{q_2}{p_2},e_0+\frac{q_3}{p_3})\cong(\Sigma
\times S^1)\cup_{(\varphi_1\cup\varphi_2\cup\varphi_3)}(V_1\cup
V_2\cup V_3).
\]

Let $\xi$ be a tight contact structure on $M$. By Theorem 1.1 of
\cite{Wu}, $\xi$ admits a vertical Legendrian circle $L$ with
twisting number $0$. We first isotope $\xi$ so that there is a
vertical Legendrian circle with twisting number $0$ in the
interior of $\Sigma \times S^1$, and each $V_i$ is a standard
neighborhood of a Legendrian circle $L_i$ isotopic to the $i$-th
singular fiber with twisting number $t_i<0$, i.e., $\partial V_i$
is convex with two dividing curves each of which has slope
$\frac{1}{t_i}$ when measured in the coordinates of $\partial V_i$
given above. Let $s_i$ be the slope of the dividing curves of
$T_i=\partial V_i$ measured in the coordinates of $T_i$. Then we
have that
\[
s_i ~=~
\left\{%
\begin{array}{ll}
    \frac{-t_iq_i+v_i}{t_ip_i-u_i} ~=~
-\frac{q_i}{p_i}+\frac{1}{p_i(t_ip_i-u_i)}, & i=1,2; \\
    \frac{-t_3(q_3+e_0p_3)+(v_3+e_0u_3)}{t_3p_i-u_3} ~=~
-e_0-\frac{q_i}{p_i}+\frac{1}{p_i(t_ip_i-u_i)}, & i=3. \\
\end{array}%
\right.
\]
We choose $t_i<<-1$ so that
$\frac{1}{b^{(i)}_0+1}<s_i<-\frac{q_i}{p_i}$ for $i=1,2$, and
$-e_0+\frac{1}{b^{(3)}_0+1}<s_3<-e_0-\frac{q_3}{p_3}$. Using the
vertical Legendrian circle $L$, we can thicken $V_i$ to $V'_i$,
s.t., $V'_i$'s are pairwise disjoint, and $T'_i=\partial V'_i$ is
a minimal convex torus with vertical dividing curves when measured
in coordinates of $T_i$. By Proposition 4.16 of \cite{H1}, there
exits a minimal convex torus $T''_i$ in the interior of $V'_i
\setminus V_i$ isotopic to $T_i$ that has dividing curves of slope
$\frac{1}{b^{(i)}_0+1}$ for $i=1,2$, and
$-e_0+\frac{1}{b^{(3)}_0+1}$ for $=3$. Let $V''_i$ be the solid
torus bounded by $T''_i$, and $\Sigma''\times S^1=M\setminus(V''_1
\cup V''_2 \cup V''_3)$.

Now we count the tight contact structures on $\Sigma''\times S^1$
and $V''_i$ that satisfy the given boundary condition. First, we
look at $V''_i$. In the coordinates in $\partial V_i$, the
dividing curves of $T''_i=\partial V''_i$ have slope
$\frac{(b^{(i)}_0+1)q_i+p_i}{(b^{(i)}_0+1)v_i+u_i}$. By part (4)
of Lemma \ref{CF} and the definitions of $u_i$, $v_i$, we have
that
$\frac{(b^{(i)}_0+1)q_i+p_i}{(b^{(i)}_0+1)v_i+u_i}=<b^{(i)}_{l_i},b^{(i)}_{l_i-1},\cdots,b^{(i)}_{2},
b^{(i)}_{1}+1>$. Thus, on each $V''_i$, there are
$|\prod^{l_i}_{j=1}(b^{(i)}_j+1)|$ tight contact structures that
satisfy the given boundary condition. Then we look at
$\Sigma''\times S^1$. The thickened torus $L_i$ bounded by
$T'_i-T''_i$ is a continued fraction block consisting of
$|b^{(i)}_0+1|$ basic slices. Let $L'_i$ be the basic slice in
$L_i$ closest to $T'_i$, and $\partial L'_i=T'_i-T'''_i$. Note
that $T'''_i$ is a minimal convex torus with dividing curves of
slope $-1$ for $i=1,2$, and $-e_0-1$ for $i=3$. Let $\Sigma'\times
S^1=M\setminus(V'_1 \cup V'_2 \cup V'_3)$. By Lemma \ref{center},
the tight contact structure on $(\Sigma'\times S^1)\cup L'_1 \cup
L'_2 \cup L'_2$ is uniquely determined by the signs of the basic
slices $L'_i$. But we can shuffle the signs of the basic slices
within a continued fraction block. Let's shuffle all the positive
signs in $L_i$ to the basic slices closest to $T'_i$. Then the
sign of $L'_i$ is uniquely determined by the number of positive
slices in $L_i$, and so is the number of positive slices in
$L_i\setminus L'_i$. Thus, the tight contact structures on
$(\Sigma'\times S^1)\cup L'_1 \cup L'_2 \cup L'_2$ and
$L_i\setminus L'_i$ are uniquely determined by these three
numbers. But there are only $|b^{(1)}_0b^{(2)}_0b^{(3)}_0|$ ways
to choose these three numbers. So there are at most
$|b^{(1)}_0b^{(2)}_0b^{(3)}_0|$ on $\Sigma''\times S^1$ that
satisfy the given boundary condition. Altogether, there are at
most $|\prod^3_{i=1}b^{(i)}_0\prod^{l_i}_{j=1}(b^{(i)}_j+1)|$
tight contact structures on $M$.

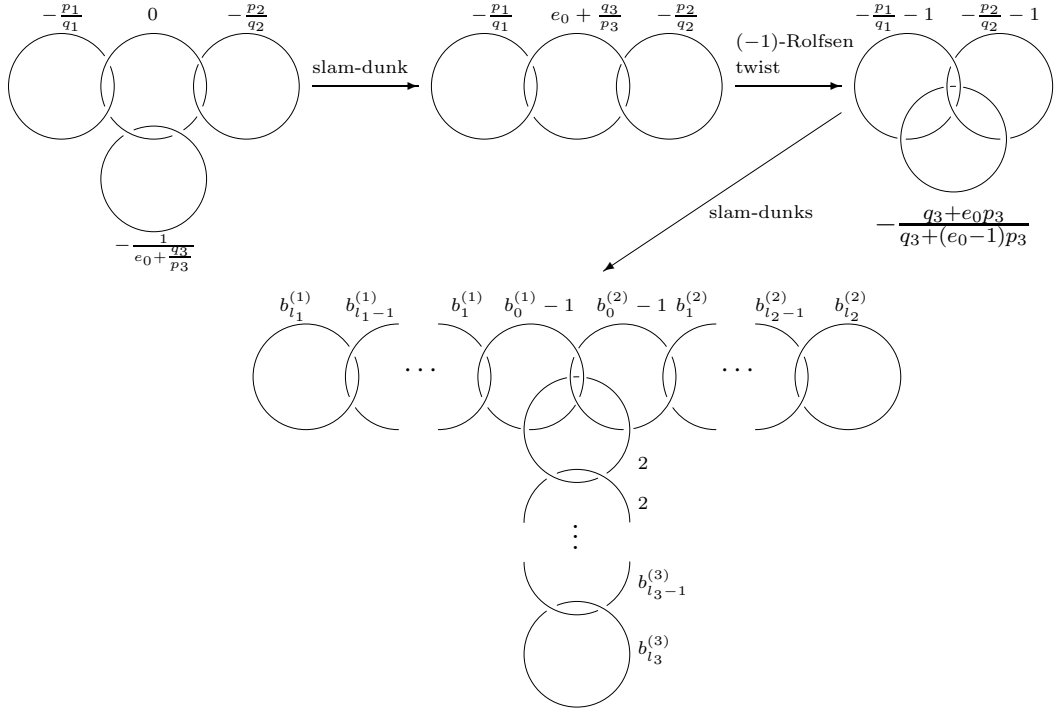
\begin{figure}

\setlength{\unitlength}{1pt}
\begin{picture}(570,290)(-210,-250)


\put(-162.5,25){\tiny{$0$}}

\put(-160,0){\arc{40}{9.04}{12.95}}

\put(-160,0){\arc{40}{6.9}{7.25}}

\put(-160,0){\arc{40}{7.50}{8.8}}

\put(-202.5,25){\tiny{$-\frac{p_1}{q_1}$}}

\put(-195,0){\arc{40}{5.9}{11.95}}

\put(-132.5,25){\tiny{$-\frac{p_2}{q_2}$}}

\put(-125,0){\arc{40}{3.72}{9.85}}

\put(-175,-62){\tiny{$-\frac{1}{e_0+\frac{q_3}{p_3}}$}}

\put(-160,-35){\arc{40}{4.33}{10.38}}


\put(-100,5){\tiny{slam-dunk}}

\put(-100,0){\vector(1,0){40}}


\put(-10,25){\tiny{$e_0+\frac{q_3}{p_3}$}}

\put(0,0){\arc{40}{9.04}{12.95}}

\put(0,0){\arc{40}{6.9}{8.8}}

\put(-40,25){\tiny{$-\frac{p_1}{q_1}$}}

\put(-35,0){\arc{40}{5.9}{11.95}}

\put(30,25){\tiny{$-\frac{p_2}{q_2}$}}

\put(35,0){\arc{40}{3.72}{9.85}}


\put(60,15){\tiny{$(-1)$-Rolfsen}}

\put(60,5){\tiny{twist}}

\put(60,0){\vector(1,0){40}}


\put(105,25){\tiny{$-\frac{p_1}{q_1}-1$}}

\put(125,0){\arc{40}{6.9}{7.88}}

\put(125,0){\arc{40}{8.05}{12.95}}

\put(145,25){\tiny{$-\frac{p_2}{q_2}-1$}}

\put(160,0){\arc{40}{3.72}{7.65}}

\put(160,0){\arc{40}{7.85}{9.8}}

\put(142.5,-20){\arc{40}{4.65}{4.75}}

\put(142.5,-20){\arc{40}{4.9}{10.8}}

\put(112.5,-55){$-\frac{q_3+e_0p_3}{q_3+(e_0-1)p_3}$}


\put(100,-10){\vector(-3,-2){90}}

\put(50,-50){\tiny{slam-dunks}}



\put(23,-145){\tiny{$2$}}

\put(0,-130){\arc{40}{4.65}{4.75}}

\put(0,-130){\arc{40}{4.9}{7.25}}

\put(0,-130){\arc{40}{7.45}{10.8}}


\put(-27.5,-85){\tiny{$b^{(1)}_0-1$}}

\put(-17.5,-110){\arc{40}{6.85}{7.88}}

\put(-17.5,-110){\arc{40}{8.05}{8.8}}

\put(-17.5,-110){\arc{40}{9.04}{12.95}}

\put(-47.5,-85){\tiny{$b^{(1)}_1$}}

\put(-52.5,-110){\arc{40}{11}{11.95}}

\put(-52.5,-110){\arc{40}{5.9}{7.85}}

\put(-65.5,-110){$\cdots$}

\put(-87.5,-85){\tiny{$b^{(1)}_{l_1-1}$}}

\put(-67.5,-110){\arc{40}{9.04}{11}}

\put(-67.5,-110){\arc{40}{7.85}{8.8}}

\put(-112.5,-85){\tiny{$b^{(1)}_{l_1}$}}

\put(-102.5,-110){\arc{40}{5.9}{11.95}}


\put(7.5,-85){\tiny{$b^{(2)}_0-1$}}

\put(17.5,-110){\arc{40}{3.72}{5.67}}

\put(17.5,-110){\arc{40}{5.9}{7.65}}

\put(17.5,-110){\arc{40}{7.85}{9.8}}

\put(37.5,-85){\tiny{$b^{(2)}_1$}}

\put(52.5,-110){\arc{40}{9.04}{11}}

\put(52.5,-110){\arc{40}{7.85}{8.8}}

\put(54.5,-110){$\cdots$}

\put(67.5,-85){\tiny{$b^{(2)}_{l_2-1}$}}

\put(67.5,-110){\arc{40}{11}{11.95}}

\put(67.5,-110){\arc{40}{5.9}{7.85}}

\put(97.5,-85){\tiny{$b^{(2)}_{l_2}$}}

\put(102.5,-110){\arc{40}{9.04}{15.1}}


\put(23,-160){\tiny{$2$}}

\put(0,-165){\arc{40}{4.33}{6.28}}

\put(0,-165){\arc{40}{9.42}{10.38}}

\put(-2,-175){$\vdots$}

\put(23,-190){\tiny{$b^{(3)}_{l_3-1}$}}

\put(0,-180){\arc{40}{7.45}{9.42}}

\put(0,-180){\arc{40}{6.28}{7.25}}

\put(23,-215){\tiny{$b^{(3)}_{l_3}$}}

\put(0,-215){\arc{40}{4.33}{10.38}}

\end{picture}

  \caption{Construction of the tight contact structures on $M$}\label{surgery2}
\end{figure}

It remains to construct
$|\prod^3_{i=1}b^{(i)}_0\prod^{l_i}_{j=1}(b^{(i)}_j+1)|$ tight
contact structures on $M$ by Legendrian surgeries of
$(S^3,\xi_{st})$. We begin with the standard surgery diagram of
$M=M(\frac{q_1}{p_1},\frac{q_2}{p_2},e_0+\frac{q_3}{p_3})$. Then,
perform a slum-dunk between the $0$-component and the
$\frac{-1}{e_0+\frac{q_3}{p_3}}$-component, after which the
$\frac{-1}{e_0+\frac{q_3}{p_3}}$-component disappears and the
original $0$-component becomes a
$(e_0+\frac{q_3}{p_3})$-component. Next we perform a
$(-1)$-Rolfson twist on the $(e_0+\frac{q_3}{p_3})$-component,
after which the three components remain trivial and have
coefficients $-\frac{p_1}{q_1}-1$, $-\frac{p_2}{q_2}-1$ and
$-\frac{q_3+e_0p_3}{q_3+(e_0-1)p_3}$. But we have
\[
-\frac{p_1}{q_1}-1=<b^{(1)}_0-1,b^{(1)}_1,\cdots,b^{(1)}_{l_1}>,
\]
\[
-\frac{p_2}{q_2}-1=<b^{(2)}_0-1,b^{(2)}_1,\cdots,b^{(2)}_{l_2}>
\]
and
\[
-\frac{q_3+e_0p_3}{q_3+(e_0-1)p_3}=<-2,\cdots,-2,b^{(3)}_0-1,b^{(3)}_1,\cdots,b^{(3)}_{l_3}>,
\]
where, on the right hand side of the last equation, there are
$e_0-1$ many $-2$'s in front of $b^{(3)}_0-1$. Now, we perform
(inverses of) the slam-dunks corresponding to these three
continued fractions here, which lead us to the diagram at the
bottom of Figure \ref{surgery2}. Note that all components in this
diagram are trivial. Since the maximal Thurston-Bennequin number
of an unknot in $(S^3,\xi_{st})$ is $-1$, it's easy to see that
there are $|\prod^3_{i=1}b^{(i)}_0\prod^{l_i}_{j=1}(b^{(i)}_j+1)|$
ways to realize this diagram by Legendrian surgeries. According to
Proposition 2.3 of \cite{Go} and Theorem 1.2 of \cite{LM}, each of
these Legendrian surgeries gives a non-isotopic tight contact
structure on $M$. Thus, up to isotopy, there are exactly
$|\prod^3_{i=1}b^{(i)}_0\prod^{l_i}_{j=1}(b^{(i)}_j+1)|$ tight
contact structures on $M$. \hfill $\Box$


\vspace{.4cm}

\begin{thebibliography}{99}
  \bibitem{CGH}
   V. Colin, E. Giroux, K. Honda,
   \emph{Notes on the isotopy finiteness,}
   arXiv:math.GT/0305210.
  \bibitem{E1}
   Y. Eliashberg,
   \emph{Classification of overtwisted contact structures ob $3$-manifolds,}
   Invent. Math. \textbf{98} (1989), no. 3, 623-637.
  \bibitem{E2}
   Y. Eliashberg,
   \emph{Contact $3$-manifolds twenty years since J. Martinet's work,}
   Ann. Inst. Fourier (Grenoble) \textbf{42} (1992), no. 1-2, 165-192.
  \bibitem{E3}
   Y. Eliashberg,
   \emph{Classification of contact structures on $\bold R\sp 3$,}
   Internat. Math. Res. Notices \textbf{1993},  no. 3, 87-91.
  \bibitem{EH}
   J. Etnyre, K. Honda,
   \emph{On the nonexistence of tight contact structures,}
   Ann. of Math. (2) \textbf{153} (2001), no. 3, 749-766.
  \bibitem{Gh}
   P. Ghiggini,
   \emph{Tight contact structures on Seifert Manifolds over $T^2$ with one singular
   fiber,}
   arXiv:math.GT/0307340.
  \bibitem{GS}
   P. Ghiggini, S. Sch\"{o}nenberger,
   \emph{On the classification of tight contact structures,}
   arXiv:math.GT/0201099.
  \bibitem{Gi1}
   E. Giroux,
   \emph{Convexit\'{e} en topologie de contact,}
   Comment. Math. Helv. \textbf{66} (1991), no. 4, 637-677.
  \bibitem{Gi2}
   E. Giroux,
   \emph{Structures de contact sur les vari\'{e}t\'{e}s fibr\'{e}es au-dessus d'une surface,}
   Comment. Math. Helv. \textbf{141} (2000), 615-689.
  \bibitem{Go}
   R. Gompf,
   \emph{Handlebody construction of Stein surfaces,}
   Ann. of Math. (2) \textbf{148} (1998), no. 2, 619-693.
  \bibitem{Ha}
   A. Hatcher,
   \emph{Notes on Basic 3-Manifold Topology,}
   http://www.math.cornell.edu/$\sim$hatcher.
  \bibitem{H1}
   K. Honda,
   \emph{On the classification of tight contact structures I,}
   Geom. Topol. \textbf{4} (2000), 309-368 (electronic).
  \bibitem{H2}
   K. Honda,
   \emph{On the classification of tight contact structures II,}
   J. Differential Geom. \textbf{55} (2000), no. 1, 83-143.
  \bibitem{HKM}
   K. Honda, W. Kazez, G. Mati\'{c},
   \emph{Convex decomposition theory,}
   Internat. Math. Res. Notices \textbf{2002},  no. 2, 55-88.
  \bibitem{K}
   Y. Kanda,
   \emph{the classification of tight contact structures on the $3$-torus,}
   Comm. Anal. Geom. \textbf{5} (1997), no. 3, 413-438.
   \bibitem{K2}
   Y. Kanda,
   \emph{On the Thurston-Bennequin invariant of Legendrian knots and non exactness of Bennequin's inequality,}
   Invent. Math. \textbf{133} (1989), 227-242.
   \bibitem{LM}
   P. Lisca, G. Mati\'c,
   \emph{Tight contact structures and Seiberg-Witten invariants,}
   Invent. Math. \textbf{129} (1997), no. 3, 509--525.
   \bibitem{Wu}
   H. Wu,
   \emph{Legendrian Vertical Circles in Small Seifert Spaces,}
   arXiv:math.GT/0310034.
\end{thebibliography}
\end{document}